\documentclass[12pt]{article}
\usepackage[cp1251]{inputenc}
\usepackage[english]{babel}
\usepackage{amsfonts}

\begin{document}

\def\C{{\bf C}}
\def\R{{\bf R}}
\def\N{{\bf N}}
\def\‘P{{\bf CP}}
\def\i{{\sqrt {-1}}}
\def\O{{\cal O}}
\def\A{{\cal A}}
\def\M{{\widehat{M}}}
\def\DW{{\widehat{\cal D}}}
\def\D{{\cal D}}
\def\Z{{\bf Z}}
\def\L{{\cal L}}
\def\A{{\cal A}}
\def\B{{\cal B}}
\def\P{{\cal P}}
\def\Pic{\rm{Pic}}
\def\Mat{\rm{Mat}}
\def\ord{\rm{ord}}
\def\Im{\rm{Im}}

\title{Commutative Rings of~Differential Operators
Connected with Two-Dimensional Abelian Varieties}
\author{A.E. Mironov
\thanks{Institute of Mathematics, 630090 Russia, Novosibirsk;
e-mail: mironov@math.nsc.ru}}
\date{}
\maketitle

\begin{abstract}
In this paper we find the explicit formulas of two dimensional
commuting ($2\times 2$)-matrix differential operators which were
introduced by Nakayashiki. The common eigen functions and
eigen values of these operators are parametrized
by the points of principally polarized Abelian varieties.
\end{abstract}

\section{Introduction}

The main result of this paper is construction of two-dimensional
($2\times 2$)-matrix differential operators
with doubly periodic coefficients from spectral functions and
the vector Baker--Akhiezer eigenfunction that was
introduced by Nakayashiki in~[1].
A~remarkable property of these operators is the fact that
they are finite-gap at every energy level. Namely,
their Bl\^och functions, the eigenfunctions of both the
differential operator and the operators of translation by periods,
are parametrized by the points of a~two-dimensional
principally polarized Abelian variety~
$X^2$.
Moreover, the eigenvalues (making the spectrum) are given by some meromorphic
function
$\lambda (z)$
on~
$X^2$
with pole on a~theta divisor. We indicate a~procedure
for constructing a~commutative ring of such operators
from an~Abelian variety~
$X^2$
(with an~irreducible theta divisor) and a~spectral function~
$\lambda$.

 The multidimensional inverse problem was first solved
by B.~A. Dubrovin, I.~M. Krichever, and S.~P. Novikov
for a~periodic Schr\"odinger operator with
a~magnetic field~[2]. They described a~procedure
for constructing such operators that are finite-gap
at one energy level, i.e., their Bl\^och functions with
a~fixed eigenvalue~
$E_0$
are parametrized by the points of a~Riemann surface
of finite genus. Later, A.~P. Veselov and
S.~P. Novikov~[3] distinguished  potential
operators (i.e., operators with zero magnetic field) among them.

The finite-gap property of the Schr\"odinger
operator at a~fixed energy level is an~exceptional phenomenon:
Feldman, Knorrer, and Trubowitz~[4] demonstrated that
the two-dimensional Schr\"odinger operator with a~smooth real potential
can be finite-gap only at one energy level.
As was first observed by Sato, for matrix operators
it is possible that the Bloch functions are parametrized
by surfaces of finite genus (in our case by the level surfaces
$\lambda = const $)
for all energy levels.
This idea was implemented by Nakayashiki~[1,\,5].

In~[1], using the Fourier--Mukai transform method~[6],
Nakayashiki constructed the Baker--Akhiezer module over a~ring
of differential operators in~
$g$
space variables,
where
$g$
is the dimension of a~principally polarized Abelian variety
$X^g$
with a~nonsingular theta divisor.
To this module there corresponds (up to conjugation)
a~commutative ring of ($g!\times g!$)-matrix
partial differential operators whose coefficients are, in general, defined locally.
Each operator corresponds to some
meromorphic function on~
$X^g$
with pole on the theta divisor
(a~spectral function) which parametrizes the eigenvalues
of the operator (these operators are now referred to as
the Nakayashiki operators). The vector-function whose components
are the elements of the basis for the Baker--Akhiezer module
parametrizes the general eigenfunctions of these operators.

Observe that this construction requires nonsingularity of
the theta divisor. As was demonstrated by Andreotti and Mayer~[7],
a~theta divisor of a~general Abelian variety
is a~nonsingular subvariety; however, a~theta divisor
of the Jacobian varieties of Riemann surfaces has singularities
for
$g>3$
and in the case of hyperelliptic surfaces also
for
$g=3$.

In Section~2 we describe the Fourier--Mukai transform~[6]
in the needed situation,
recall Naka\-yashi\-ki's construction of the Baker--Akhiezer module~[1], and
indicate connection between Krichever's construction~[8]
and the Fourier--Mukai transform~[6].

In Section~3 we give a~new short analytical
proof of Nakayashiki's theorem of freeness of
the Baker--Akhiezer module
in dimension~2.
It  essentially uses the fact that the Abelian variety
$X^2$
in question
is two-dimensional unlike
the proof of the general case in~[1] and requires a~minimal
apparatus of algebraic geometry.
Here we also introduce a~basis for the Baker--Akhiezer module for
$g=2$
in which we manage to find the coefficients of the Nakayashiki operators.

In Section~4 we describe an~effective procedure for constructing
the Nakayashiki operators for
$g=2$.
In Propositions 1--3 we obtain
explicit formulas for operators generating, as proven in
Lemma~9, the whole ring of the Nakayashiki operators.
In Proposition~1 we
introduce commuting operators
$Z_1,\dots,Z_g$
such that, for every Nakayashiki operator~
$L$,
the commutator
$[L,Z_j]$
is a~Nakayashiki operator too. In Proposition~2
we find explicit formulas for the second-order Nakayashiki operators
and in Proposition~3, explicit formulas for the operators~
$Z_j$.

The author is grateful to S.~P. Novikov
and I.~A. ~Taimanov for posing the problem
and to I.~A. Taimanov for useful discussions and remarks.

\section{Nakayashiki's Construction}

Let
$X^g={\C}^g/({\Z}^g+\Omega {\Z}^g)$
be a~principally polarized  Abelian complex variety, where
$\Omega$
is a~symmetric ($g\times g$)-matrix
with
$\Im \Omega>0$.
Denote by
$\Pic ^0 (X^g)$
the Picard variety of~
$X^g$.  If
$X^g$
is principally polarized then
$\Pic ^0 (X^g)$
is isomorphic to~
$X^g$.
The isomorphism is given by the mapping
$z\rightarrow x$,
where
$z\in X^g$ and $x \in \Pic ^0(X^g)$.
The sections of a~bundle
$x$,
lifted to the universal covering
${\C}^g$,
are given by functions
$h(z)$
satisfying the periodicity conditions
$$
h(z+\Omega m+n)=\exp(-2\pi i \langle m,x \rangle )h(z),
$$
where
$m,n \in{\Z}^g$
and
$\langle m,x \rangle=m_1x_1\dots+m_gx_g$.

Denote by
${\cal P}$
the Poincar\'e bundle over
$X^g\times \Pic ^0(X^g)$.
It is determined by the following properties. The bundle corresponding to
$x\in \Pic ^0(X^g)$
is isomorphic to
${\cal P}\vert_{X^g\times\{x\}}$
and the bundles
${\cal P}\vert_{X^g\times\{0\}}$
and
${\cal P}\vert_{\{0\}\times \Pic ^0(X^g)}$
are trivial. The sections of
${\cal P}$,
lifted to the universal covering
${\C}^g\times{\C}^g$,
are given by functions
$f(z,x)$
such that
$$
f(z+\Omega m_1+n_1,x+\Omega
m_2+n_2)= \exp(-2\pi i(\langle m_1,x\rangle + \langle
m_2,z\rangle))f(z,x),
                            \eqno{(1)}
$$
where
$m_j, n_j\in {\Z}^g$.
Denote by
$z=(z_1,\dots,z_g)^{\top}\in{\C}^g$
the coordinates of points on the universal covering of~
$X^g$
and by
$x=(x_1,\dots,x_g)^{\top}\in{\C}^g$,
the coordinates of points on the universal covering of~
$\Pic ^0(X^g)$.
We identify the sections of
${\cal P}$
with the functions on
${\C}^g\times{\C}^g$
satisfying ~(1).

The theta function with characteristic
$[a, b]$
is defined by the series
$$
\theta[a,b](z,\Omega)=\sum\limits_{n\in {\Z}^g}\exp
(\pi i \langle \Omega (n+a),(n+a) \rangle
+2\pi i
\langle (n+a),(z+b)\rangle ),
$$
where
$a, b\in {\C}^g$.
For short,
the function
$\theta(z)=\theta[0,0](z,\Omega)$
is referred to as a~{\it theta function}.
The function
$\theta[a, b](z,\Omega)$
possesses the periodicity properties
$$
\theta[a,b](z+m,\Omega)=\exp(2\pi i
\langle a,m\rangle)
\theta[a,b](z,\Omega),
$$
$$ \theta[a,b](z+\Omega m,
\Omega)=\exp(-2\pi i
\langle b,m \rangle
-\pi i \langle m\Omega ,m\rangle -2\pi i
\langle m,z \rangle)\theta[a,b](z,\Omega),
$$
$m\in {\Z}^g.$

Denote by
$\Theta$
the zero set of a~theta function
(a~theta divisor)
which is a~subvariety in~
$X^g$.

Given
a~subvariety $Y\subset X^g$ of codimension~1, we let
$A_Y$
stand for  the space of meromorphic functions on~
$X^g$
with pole in~
$Y$.
We denote by
$F_Y(U)$
the space of meromorphic sections of
${\cal P}$
over subsets of the form
$X^g\times U\subset X^g\times \Pic ^0(X^g)$
which have pole in
$Y\times U$,
where
$U$
is an~open subset in~
$\Pic ^0(X^g)$.

The set
$F_Y=\bigcup\nolimits_U F_Y(U)$
is called the Fourier--Mukai transform of
$A_Y$~[6].

In~[1] Nakayashiki constructed ``covariant differentiation''
operators (a~connection on~
${\cal P}$)
$$
\nabla_j: F_Y(U)\rightarrow F_Y(U),
\quad
\nabla_k\nabla_j=\nabla_j\nabla_k,
\quad
k,j=1,\dots,g,
$$
that give
$F_Y(U)$
a~ module structure over the ring
${\cal O}_U[\nabla_1,\dots,\nabla_g]$,
where
$U\subset \Pic ^0(X^g)$
is an~open subset and
${\cal O}_U$
is the ring of analytic functions on~
$U$.
By construction,
$F_Y(U)$
is also an~$A_Y$-module.

Nakayashiki introduced the Baker--Akhiezer functions
$$
f(z,x)
\exp\Biggl(-\sum\limits_{j=1}^gx_j\partial_{z_j}\log\theta(z)\Biggr),
\quad  f(z,x)\in F_{\Theta}.
$$
The operators
$\nabla_j$
are defined by the formula
$$
\nabla_j=\partial_{x_j}-\partial_{z_j}\log\theta(z).
$$

Let
$M_c$
denote the space of the Baker--Akhiezer functions
$$
\Biggl\{ f(z,c+x)\exp
\Biggl(-\sum\limits_{j=1}^gx_j\partial_{z_j}\log
\theta(z)\Biggr)\mid f(z,c+x)\in \bigcup_U F_{\Theta}(c+U)\Biggr\},
$$
where
$U\subset {\C}^g$
is a~neighborhood of zero,
$c\in{\C}^g$,
and
$f(z,c+x)$
belongs to the range of the Fourier--Mukai transform~
$F_{\Theta}$.
Each function
$f(z,c+x)$
is defined on
${\C}^g\times U$
and has pole only in
$\Theta\times U$.
Roughly speaking,
$f(z,c+x)$
is a~germ of a~function at
$x=0$
with respect to~
$x$.
It follows from (1) that
$$
f(z+\Omega m +n,c+x)=\exp(-2\pi i\langle m,(c+x)\rangle )f(z,c+x).
$$
We denote by
$M_c(k)$
the subset of functions in~
$M_c$
such that
$f(z,c+x)$
has a~pole of order
$\leq k$
in
$\Theta\times U$.
The definition implies that
$M_c$
is a~module over the ring
$
{\cal D}={\cal O}[\partial_{x_1},\dots,\partial_{x_g}]
$
of differential operators
(a~${\cal D}$-module),
where
${\cal O}$
is the ring of analytic functions in
$x_1,\dots,x_g$
in a~neighborhood of
$0\in{\C}^g$. The ${\cal D}$-module
$M_c$
is referred to as the {\it Baker--Akhiezer module}.
From the definition of the Fourier--Mukai transform we infer that
$M_c$
also possesses the structure of an~$A_{\Theta}$-module.
 The ${\cal D}$-module $M_c$ can be described
by means of theta functions~[1].

{\bf Lemma 1.}
{\sl
The equality holds:
$$
M_c=\sum\limits_{n=1}^{\infty}\sum\limits_{a\in {\Z}^g/n{\Z}^g}
{\cal O} \cdot
\frac{\theta[\frac{a}{n},0](nz+c+x,n\Omega)}{\theta^n(z)}
\exp\Biggl(-\sum\limits_{k=1}^gx_k
\partial_{z_k}\log\theta(z)\Biggr).
$$
}

We also need the following ~[1]:

{\bf Lemma 2.}
{\sl
The identity holds:
$$
\partial_{x_j}\Biggl(
\frac{\theta[\frac{a}{n},0](nz+c+x,n\Omega)}{\theta^n(z)}
\exp\Biggl(-\sum\limits_{k=1}^gx_k
\partial_{z_k}\log\theta(z)\Biggr)\Biggr)
$$
$$
=\frac{1}{n}
\partial_{z_j}\left(
\frac{\theta[\frac{a}{n},0](nz+c+x,n\Omega)}{\theta^n(z)}
\right)
\exp\Biggl(-\sum\limits_{k=1}^gx_k\partial_{z_k}\log\theta(z)\Biggr),
$$
$j=1,\dots,g.$
}

 The next assertion was proven in~[1].

{\bf Nakayashiki's Theorem.}
{\sl
If
$\Theta$
is a~nonsingular variety and
$c\ne 0$
then $M_c$
is a~free
${\cal D}$-module
of rank~
$g!$.
}

Fix a~basis
$\Phi_c=(\phi_{1c}(z,x),\dots,\phi_{g!c}(z,x))^{\top}$
for the ${\cal D}$-module $ M_c$.
Suppose that
$\lambda (z)\in A_\Theta$.
Since
$M_c$
is a~free ${\cal D}$-module,
there exist unique operators
$[L_{\Phi_c}(\lambda )]_{kj}\in{\cal D}$
such that
$$
\sum\limits_{j=1}^{g!}[L_{\Phi_c}(\lambda )]_{kj}\phi_{jc}=\lambda \phi_{kc},
\quad
k=1,\dots,g!.
$$
Consequently,
$$
L_{\Phi_c}(\lambda)\Phi_c=\lambda\Phi_c,
                            \eqno{(2)}
$$
with
$[L_{\Phi_c}(\lambda)]_{kj}$
the entries of the matrix differential operator
$L_{\Phi_c}(\lambda)$
and
$\lambda\Phi_c=(\lambda\phi_{1c},\dots,\lambda\phi_{g!c})^{\top}$.
Since
$L_{\Phi_c}(\lambda)$
are differential operators in
$x_j$
while $\lambda$
depends only on
$z$,
from (2) we obtain the commutation condition
$$
L_{\Phi_c}(\lambda\mu)=L_{\Phi_c}(\lambda)L_{\Phi_c}(\mu)=
L_{\Phi_c}(\mu)L_{\Phi_c}(\lambda),
$$
where
$\mu(z)\in A_{\Theta}$.
We arrive at the following

{\bf Corollary \rm[1].}
{\sl
There is a~ring embedding
$$
L_{\Phi_c}:A_\Theta\rightarrow \Mat (g!,{\cal D}),
$$
where
$\Mat (g!,{\cal D})$
is the ring of ($g! \times g!$)-matrix
differential operators. The image of the embedding is
a~commutative ring of differential operators.
}

As indicated in~[5], Nakayashiki's construction
generalizes Krichever's construction~[8] which
can be interpreted in terms of the Fourier--Mukai transform
as follows:

Recall the construction of the Baker--Akhiezer function~[8]. Let
$\Gamma$
be a~Riemann surface of genus~
$g$, let
$D=p_1+\dots +p_g$
be a~nonspecial positive divisor on~
$\Gamma$,
and  let
$\infty$
be a~point on
$\Gamma$
other than the points of the divisor.  Take a~local parameter
$k^{-1}$ at $\infty$
so that
$k^{-1}(\infty)=0$.
The {\it one-point Baker--Akhiezer function with spectral data}
$\{\Gamma,\infty,p_1,\dots,p_g,k^{-1}\} $
is a~function
$\psi(z,x)$, $z\in\Gamma$,
defined to within a~factor
depending only on~$x$, by the following properties:

(1) $\psi(z,x)$
is meromorphic on
$\Gamma\backslash \infty$
and the poles of $\psi$ are independent of~$x$ and
coincide with
$\{p_1,\dots, p_g\}$;

(2) the function
$\psi(z,x)\exp(-kx)$
is analytic in a~neighborhood of~
$\infty$.

This function has the form
$$
\psi(z,x)=
\frac{\theta\Bigl({\cal A}(z)-\sum\limits_{j=1}^g{\cal
A}(p_j)-\Delta+Vx\Bigr)} {\theta\Bigl({\cal
A}(z)-\sum\limits_{j=1}^g{\cal A}(p_j)-\Delta\Bigr)}
\exp\biggl(2\pi ix\int\limits_{p_0}^{z}\eta\biggr),
$$
where
${\cal A}: \Gamma \to X^g$
is the Abel mapping with a~basepoint
$p_0$, $X^g$
is the Jacobian variety of
$\Gamma$, $\Delta$
is the vector of Riemann constants, and
$V\in{\C}^g$
is a~vector determined from the spectral data.

The function
$$
\frac{\theta\Bigl({\cal A}(z)-\sum\limits_{j=1}^g{\cal
A}(p_j)-\Delta+Vx\Bigr)} {\theta\Bigl({\cal
A}(z)-\sum\limits_{j=1}^g{\cal A}(p_j)-\Delta\Bigr)}
$$
is the
restriction to~
${\cal A}(\Gamma)\times\{Vx\}$ of the function
$$
\frac{\theta\Bigl(\tilde{z}-\sum\limits_{j=1}^g{\cal
A}(p_j)-\Delta+ \tilde{x}\Bigr)}
{\theta\Bigl(\tilde{z}-\sum\limits_{j=1}^g{\cal A}(p_j) -
\Delta\Bigr)}
$$
which belongs  to the range of the
Fourier--Mukai transform
$F_{\Theta'}$,
where
$\Theta'\subset
X^g
$
is the subvariety defined by the equation
$$
\theta\biggl(\tilde{z}-\sum\limits_{j=1}^g{\cal
A}(p_j)-\Delta\biggr)=0.
$$

Denote by
$\widetilde{{\cal D}}(\psi(z,x))$
the
$\widetilde{{\cal D}}$-module
$\{d\psi(z,x)\mid  d\in \widetilde{{\cal D}}\}$,
where
$\widetilde{{\cal D}}$
is the ring of differential operators in~
$x$.
As shown in~[8],
$\widetilde{{\cal D}}(\psi(z,x))$
is a~free
$\widetilde{{\cal D}}$-module; moreover,
for every meromorphic function~
$f(z)$
on~
$\Gamma$
with a~sole pole at~
$\infty$,
there is a~unique differential operator~
$L(f)$
such that
$$
L(f)\psi(z,x)=f(z)\psi(z,x).
$$
Hence, we obtain a~correspondence between the
spectral data of commutative rings of
scalar differential operators in the
Burchnall--Chaundy--Krichever theory
and the spectral data of commutative rings of the Nakayashiki matrix
differential operators:
$$
\{\Gamma,\infty,D,f\}
\longleftrightarrow
\{X^g,\Theta,c,\lambda\}.
$$

\section{Proof of Nakayashiki's Theorem  for $g=2$}

Here we give a~new proof of Nakayashiki's theorem for
$g=2$.
Introduce the following functions in~
$M_c$:
$$
\psi=\frac{\theta(z+c+x)}{\theta(z)}
\exp(-x_1\partial_{z_1}\log\theta(z)
-x_2\partial_{z_2}\log\theta(z)),
$$
$$
\psi_{c'}=\frac{\theta(z+c+c'+x)\theta(z-c')}{\theta^2(z)}
\exp(-x_1\partial_{z_1}\log\theta(z)
-x_2\partial_{z_2}\log\theta(z)),
$$
and define the surface
$$
\Gamma_{c'}=\{z\in{\C}^2|\ \theta(z-c')=0\}.
$$
In Lemmas~4 and 5 we prove that
$
{\cal D}(\psi, \psi_{c'})=\{d_1\psi+d_2\psi_{c'}\mid  d_1,
d_2\in{\cal D}\} $ is a~free ${\cal D}$-module.
In Lemma~6 we
show that the
${\cal D}$-modules
${\cal D}(\psi, \psi_{c'})$ and
$M_c$
coincide.

Henceforth subscripts of a~theta function denote differentiation
with respect to the corresponding variable:
$\theta_j(z)=\partial_{z_j}\theta(z)$
and
$\theta_{kj}(z)=\partial_{z_k}\partial_{z_j}\theta(z).$

{\bf Lemma 3.}
{\sl
The functions
$a\theta_1(z)+b\theta_2(z)$
and
$\theta(z+e)$,
where
$a, b\in{\C}$ and $e\in{\C}^2$,
have two zeros on
$\Theta$
with multiplicity counted to within elements of the lattice
${\Z}^2+\Omega {\Z}^2$.
}

{\sc Proof.}
Let
$X^2$
be the Jacobian variety of some Riemann surface~
$\Gamma$
of genus~2. Denote by
$\widehat{\Gamma}$
the universal covering of~
$\Gamma$.
The fundamental group
$\pi_1(\Gamma)$
acts on~
$\widehat{\Gamma}$ and we
denote
a~fundamental domain by ~$S$.
Let
${\cal A}:\widehat{\Gamma}\rightarrow X^2$
be the Abel mapping. The range
${\cal A}(\widehat{\Gamma})$
has the equation
$\theta(z-\Delta)=0.$
It is easy to verify that the sought number of zeros equals
$$
\frac{1}{2\pi i}\int\limits_{\partial S}
d\log(a\theta_1({\cal A}(p)-\Delta)+b\theta_2({\cal
A}(p)-\Delta))= \frac{1}{2\pi i}\int\limits_{\partial S}
d\log\theta({\cal A}(p)+e-\Delta)=2.
$$
Lemma 3 is proven.

{\bf Lemma 4.}
{\sl
If
$c'\in{\C}^2$
and
$c'\not\in{\Z}^2+\Omega {\Z}^2$
then
$\psi_{c'}$
does not belong to the set of functions
$$
\{ \alpha_1(x)\partial_{x_1}\psi+
\alpha_2(x)\partial_{x_2}\psi+ \alpha_3(x)\psi\mid
 \alpha_1(x),\alpha_2(x),\alpha_3(x)\in {\cal O}\}.
$$
}

{\sc Proof.}
Suppose that the assertion of the lemma is false.
Then there exist
$\alpha_1(x),\alpha_2(x),\alpha_3(x)\in{\cal O}$
such that
$$
\alpha_1(x)\partial_{x_1}\psi+
\alpha_2(x)\partial_{x_2}\psi+\alpha_3(x)\psi=\psi_{c'}.
$$
Thereby the following equality holds for
$z\in\Gamma_{c'}$:
$$
\alpha_1(x)\left( \frac{\theta_1(z+c+x)} {\theta(z)}-
\frac{\theta(z+c+x)\theta_1(z)}
{\theta^2(z)}\right)
$$
$$
+\alpha_2(x)\left(
\frac{\theta_2(z+c+x)}
{\theta(z)}-
\frac{\theta(z+c+x)\theta_2(z)}
{\theta^2(z)}\right)+
\alpha_3(x)
\frac{\theta(z+c+x)}
{\theta(z)}=0.
$$
Let
$p$
be an~intersection point of
$\Gamma_{c'}$
and
$\Theta$
(Lemma~3).
It follows from the last equality that
$$
\alpha_1(x)\theta_1(p)+\alpha_2(x)\theta_2(p)=0.
$$

Since
$\Theta$
is a~smooth Riemann surface, we have
$\theta_1(p)\ne 0$
or
$\theta_2(p)\ne 0$.
For definiteness, assume that
$\theta_1(p)\ne 0$.
Then
$$
\frac{\alpha_1(x)}{\alpha_2(x)}=-\frac{\theta_2(p)}{\theta_1(p)}.
$$
Consequently,
$$
-\frac{\theta_2(p)}{\theta_1(p)}
\partial_{x_1}\log\theta(z+c+x)+\partial_{x_2}\log\theta(z+c+x)
=-\frac{\theta_2(p)}{\theta_1(p)}
\frac{\theta_1(z)}{\theta(z)}+
\frac{\theta_2(z)}{\theta(z)}-
\frac{\alpha_3(x)}{\alpha_2(x)}.
$$
The pole (in ~
$x$)
of the left-hand side depends on
$z\in\Theta$,
whereas the pole of the right-hand side does not; a~contradiction.
Lemma~4 is proven.

{\bf Lemma 5.}
{\sl
${\cal D}(\psi, \psi_{c'})$
is a~free
${\cal D}$-module
of rank~$2$.
}

{\sc Proof.}
Suppose that
${\cal D}(\psi, \psi_{c'})$
is not a~free
${\cal D}$-module.
Then there exist operators
$d_1, d_2\in{\cal D}$
such that
$$
d_1\psi_{c'}+d_2\psi=0.
                            \eqno{(3)}
$$
Consider the case in which
$\ord (d_1)>\ord (d_2)-1$, where $\ord $
is the order of the operator. Suppose that
$\ord (d_1)=n$.
The operator
$d_1$
looks like
$$
d_1=f_n(x)\partial^n_{x_1}+
f_{n-1}(x)\partial_{x_1}^{n-1}\partial x_2+
\dots+
f_0(x)\partial_{x_2}^n+\dots,
$$
where
$f_j(x)\in{\cal O},\ j=0,\dots,n$.
Divide ~(3) by
$$
\exp(-x_1\partial_{z_1}\log\theta(z)-
x_2\partial_{z_2}\log\theta(z)),
$$
multiply by
$\theta^{n+2}(z)$,
and assume
$z\in\Theta$.
We obtain
$$
\theta(z+c+c'+x)\theta(z-c')
(f_n(x)\theta_1^n(z)+f_{n-1}(x)\theta_1^{n-1}(z)\theta_2(z)+\dots+
f_0(x)\theta_2^n(z))=0.                     \eqno(4)
$$
By Lemma~3, there is a~point
$p\in{\C}^2$
such that
$p\in\Theta$
and
$\theta_1(p)=0$.
Put $z=p$ in~(4).
We obtain
$f_0=0$.
Divide (4) by
$\theta_1(z)$
and put again
$z=p$.
We obtain
$f_1=0$.
Proceeding similarly, we obtain
$f_n=f_{n-1}=\dots =f_0=0$.
Consequently, the inequality
$\ord (d_1)>\ord (d_2)-1$
is impossible. Similarly, we show that the inequality
$\ord (d_1)+1<\ord (d_2)$
is impossible either. Consider the case in which
$\ord (d_1)+1=\ord (d_2)=n$.
Let the operators
$d_1$
and
$d_2$
look like
$$
d_1=f_{n-1}(x)\partial_{x_1}^{n-1}+
f_{n-2}(x)\partial_{x_1}^{n-2}\partial_{x_2}+
\dots+
f_0(x)\partial_{x_2}^{n-1}+\dots,
$$
$$
d_2=g_n(x)\partial_{x_1}^n+
g_{n-1}(x)\partial_{x_1}^{n-1}\partial_{x_2}+
\dots+
g_0(x)\partial_{x_2}^n+\dots.
$$
Divide ~(3) by
$$
\exp(-x_1\partial_{z_1}\log\theta(z)-
x_2\partial_{z_2}\log\theta(z)),
$$
multiply by
$\theta^{n+1}(z)$,
and assume
$z\in\Theta$.
Then
$$
\theta(z+c+c'+x)\theta(z-c')
(f_{n-1}(x)\theta_1^{n-1}(z)+
f_{n-2}(x)\theta_1^{n-2}(z)\theta_2(z)+\dots+
f_0(x)\theta_2^{n-1}(z))
$$
$$
-\theta(z+c+x)(g_n(x)\theta_1^n(z)+
g_{n-1}(x)\theta_1^{n-1}(z)\theta_2(z)+\dots+
g_0(x)\theta_2^n(z))=0.
$$
Decompose
$$
f_{n-1}(x)\theta_1^{n-1}(z)+
f_{n-2}(x)\theta_1^{n-2}(z)\theta_2(z)+\dots+
f_0(x)\theta_2^{n-1}(z),
$$
$$
g_n(x)\theta_1^n(z)+
g_{n-1}(x)\theta_1^{n-1}(z)\theta_2(z)+\dots+
g_0(x)\theta_2^n(z)
$$
into factors to obtain
$$
\theta(z+c+c'+x)\theta(z-c')
(a_{n-1}(x)\theta_1(z)+b_{n-1}(x)\theta_2(z))\dots
(a_1(x)\theta_1(z)+b_1(x)\theta_2(z))
$$
$$
-\theta(z+c+x)
(a'_n(x)\theta_1(z)+b'_n(x)\theta_2(z))\dots
(a'_1(x)\theta_1(z)+b'_1(x)\theta_2(z))=0,
\eqno{(5)}
$$
where
$a'_j(x)$, $b'_j(x)$, $a_k(x)$, and $b_k(x)$
are some functions,
$j=1,\dots,n$, $k=1,\dots,n-1$.

Fix a~point
$x$.
Observe that if two functions
$a'_j\theta_1(z)+b'_j\theta_2(z)$
and
$a_k\theta_1(z)+b_k\theta_2(z)$
have a~common zero on~
$\Theta$
then they are proportional. Indeed, suppose the contrary.
Then the ratio of these functions  has a~single simple pole on
$\Theta\subset X^2$
(by Lemma~3). Consequently, the mapping
$$
\bigl(a'_j\theta_1(z)+b'_j\theta_2(z): a_k\theta_1(z)+b_k\theta_2(z)\bigr):
\Theta\rightarrow {\bf CP}^1,
$$
where
${\bf CP}^1$
is the one-dimensional projective space, is
an~isomorphism. Since the genus of the Riemann surface
$\Theta\subset X^2$
equals~2, this is a~contradiction.

The zeros of
$a_k\theta_1(z)+b_k\theta_2(z)$
and
$\theta(z+c+x)$
on
$\Theta$
cannot coincide, since these functions
can be regarded as sections of the line bundle over
$\Theta\subset X^2$
and their ratio
$$
\frac{a_k\theta_1(z)+b_k\theta_2(z)}{\theta(z+c+x)}
                                \eqno{(6)}
$$
is again a~section of the line bundle but it has neither
poles nor zeros. The only sections with such properties are constants,
but ~(6) is not constant on~
$\Theta$.

Consequently, (5) amounts to the equality
$$
\theta(z+c+c'+x)\theta(z-c')+
\theta(z+c+x)(a(x)\theta_1(z)+b(x)\theta_2(z))=0,
$$
where
$a(x)$ and $b(x)$
are some functions. We conclude that the functions
$\theta (z-c')$
and
$\theta (z+c+x)$
must have a~common zero on
$\Theta$
for every
$x$,
which is impossible as we easily see, for instance,
from Riemann's theorem about the zeros of a~theta function
on a~Riemann surface. Lemma~5 is proven.

To complete the proof of Nakayashiki's theorem, it suffices
to demonstrate that the ${\cal D}$-modules
${\cal D}(\psi, \psi_{c'})$
and
$M_c$
coincide.

{\bf Lemma 6.}
{\sl The following equality holds:
${\cal D}(\psi, \psi_{c'})=M_c$.
}

{\sc Proof.}
The inclusion
${\cal D}(\psi, \psi_{c'})\subset M_c$
is obvious. Fix a~point
$x$.
Denote by
${\cal D}^k(\psi, \psi_{c'})\subset {\cal D}(\psi, \psi_{c'})$
the set of functions
$$
\{d_1\psi_{c'}+d_2\psi\mid  d_1, d_2\in{\cal D},\ \ord
(d_1)\leq k-2,\ \ord (d_2)\leq k-1\}.
$$
To prove the lemma,
it suffices to demonstrate that
$$
\dim_{{\C}}M_c(k)=
\dim_{{\C}}{\cal D}^k(\psi, \psi_{c'}).
$$
By Lemma~1,  we have
$$
\dim_{{\C}}M_c(k)=
\dim_{{\C}}
\biggl\{\sum\limits_{a\in {\Z}^g/k{\Z}^g}
\beta_{a}
\theta\biggl[\frac{a}{k},0\biggr](kz+c+x,k\Omega)\mid  \beta_{a}\in{\C}\biggr\}.
$$
Consequently,
$\dim_{{\C}}M_c(k)=k^2$
(see, for instance, [9]). Since the dimension of the space
of differential operators  of order
$\leq k$ with constant coefficients
equals
$$
1+\dots+(k+1)=\frac{(k+1)(k+2)}{2}
$$
and the ${\cal D}$-module
${\cal D}(\psi, \psi_{c'})$
is free, we have
$$
\dim_{{\C}}{\cal D}^k(\psi, \psi_{c'})=
\frac{(k-1)k}{2}+\frac{k(k+1)}{2}=k^2.
$$
Lemma 6 is proven.

\section{A~Commutative Ring of ($2 \times 2$)-Matrix
Differential Operators}

Denote by
$D$
the ring of differential operators
$
{\C}[\partial_{z_1},\dots, \partial_{z_g}]
$
with constant coefficients. By definition,
$M_c$
is a~$D$-module.
Let
$\Phi_c=(\phi_{1c},\dots,\phi_{g!c})^{\top}$
be a~basis for~
$M_c$.
Then we have the ring embedding
$$
{\cal L}_{\Phi_c}:  D\rightarrow   \Mat (g!,{\cal D})
$$
defined by the formula
$$
{\cal L}_{\Phi_c}(d)\Phi_c=d\Phi_c,
$$
where
$ d\in D$ and $d\Phi_c=(d\phi_{1c},\dots,d\phi_{g!c})^{\top}$.
The range of
${\cal L}_{\Phi_c}$
is a~commutative ring of matrix differential operators which is isomorphic to~
$D$.
Denote
the operator
${\cal L}_{\Phi_c}\left( \partial_{z_j} \right)$ by $Z_j$,
$j=1,\dots,g.$

Find the commutator
$[L_{\Phi_c}(\lambda),Z_j]$
for
$\lambda\in A_{\Theta}$:
$$
Z_jL_{\Phi_c}(\lambda)\Phi_c=Z_j(\lambda\Phi_c)=
\lambda(Z_j\Phi_c)=\lambda(\partial_{z_j}\Phi_c),
$$
$$
L_{\Phi_c}(\lambda)Z_j\Phi_c=L_{\Phi_c}\partial_{z_j}\Phi_c=
\partial_{z_j}(L_{\Phi_c}(\lambda)\Phi_c)=
\partial_{z_j}(\lambda\Phi_c)=
(\partial_{z_j}\lambda)\Phi_c+
\lambda(\partial_{z_j}\Phi_c).
$$
Thus,
$$
[L_{\Phi_c}(\lambda),Z_j]\Phi_c=(\partial_{z_j}\lambda)\Phi_c,
$$
which proves the following

{\bf Proposition 1.}
{\sl
$[L_{\Phi_c}(\lambda),Z_j]=L_{\Phi_c}(\partial_{z_j}\lambda)$
for
$\lambda\in A_{\Theta}$, $j=1,\dots,g$.
}

We now clarify how the matrix operators change when we replace a~
basis for the ${{\cal D}}$-module~
$M_c$.
Suppose that
$\Phi_c=(\phi_{1c},\dots,\phi_{g!c})^{\top}$
and
$\Psi_c=(\psi_{1c},\dots,\psi_{g!c})^{\top}$
are two bases for the ${{\cal D}}$-module~
$M_c$.
By Nakayashiki's theorem, there exist unique operators
${\cal A}, {\cal B}\in \Mat (g!,{\cal D})$
such that
${\cal A}\Phi_c=\Psi_c$ and ${\cal B}\Psi_c=\Phi_c$;
hence,
${\cal A}{\cal B}={\cal B}{\cal A}$
is the identity in~
$\Mat (g!,{\cal D})$.
We denote
the operator
${\cal B}$ by ${\cal A}^{-1}$.
It is easy to verify that
$$
{\cal A} L_{\Phi_c}(\lambda){\cal A}^{-1}\Psi_c =\lambda\Psi_c,
\quad
{\cal A} {\cal L}_{\Phi_c}(d){\cal A}^{-1}\Psi_c =d\Psi_c,
$$
where
$\lambda\in A_{\Theta}$ and $d\in D$.
Consequently, we have the following

{\bf Lemma 7}
{\sl
The equalities hold:
$$
{\cal A} L_{\Phi_c}(\lambda){\cal A}^{-1}=L_{\Psi_c}(\lambda),
 \quad
{\cal A} {\cal L}_{\Phi_c}(d){\cal A}^{-1}={\cal L}_{\Psi_c}(d).
$$
}

Henceforth we suppose that
$g=2$.
In Proposition~2 we find the second-order operators~
$L_{\Phi_c}(\lambda)$.
They correspond to the following spectral functions:
$$
\lambda=\frac{\partial^2 \log\theta(z)}
{\partial z_k\partial z_j}, \quad  j,k=1,2.
$$

In Proposition~3 we find the operators
$Z_1$ and~$Z_2$.
In Lemma~9 we prove that the operators of Propositions 1--3
generate the ring
$L_{\Phi_c}(A_{\Theta})$.

The following lemma can be proven by straightforward calculations
which we omit:

{\bf Lemma 8.}
{\sl
The  equality holds:
$$
\partial_{x_k}\partial_{x_j}
\left(\frac{\theta(z+c+x)}{\theta(z)}
\exp(-x_1\partial_{z_1}\log\theta(z)
-x_2\partial_{z_2}\log\theta(z))\right)
$$
$$
=\partial_{z_k}\partial_{z_j}
\left(\frac{\theta(z+c+x)}{\theta(z)}\right)
\exp(-x_1\partial_{z_1}\log\theta(z)-x_2
\partial_{z_2}\log\theta(z))
$$
$$
+\partial_{z_k}\partial_{z_j}
(\log\theta(z))\frac{\theta(z+c+x)}{\theta(z)}
\exp(-x_1\partial_{z_1}\log\theta(z)-x_2
\partial_{z_2}\log\theta(z)),
$$
where
$k,j=1,2$.
}

Let
$c'$
be in a~general position.
Denote by
$p_1, p_2\in{\C}^2$
the  intersection points of~
$\Gamma_{c'}$
with
$\Theta$
($p_1$
and
$p_2$
are defined to within elements of the lattice
${\Z}^2+\Omega {\Z}^2 )$.
Denote by
$L_{c,c'}$
the Nakayashiki operators for the basis
$\psi, \psi_{c'}$.
For brevity, we introduce the notations
$$
H_{c,c'}^{kj}=
[L_{c,c'}(\partial_{z_k}\partial_{z_j}\log\theta(z))]_{11},
\quad
F^{kj}_{c,c'}=
[L_{c,c'}(\partial_{z_k}\partial_{z_j}\log\theta(z))]_{12}.
$$

{\bf Proposition 2.}
{\sl
The  equality holds:
$$
[L_{c,c'}(
\partial_{z_k}\partial_{z_j}\log\theta(z))]_{11}=
-\partial_{x_k}\partial_{x_j}
+f_{c,c'}^{kj}\partial_{x_1}+g_{c,c'}^{kj}
\partial_{x_2}+h_{c,c'}^{kj},
$$
with
$$
f_{c,c'}^{kj}=\frac{\theta_2(p_1)\theta_2(p_2)}
{\theta_2(p_2)\theta_1(p_1)-\theta_2(p_1)\theta_1(p_2)}
$$
$$
\times\left(
\partial_{x_k}\log\theta(p_1+c+x)\frac{\theta_j(p_1)}{\theta_2(p_1)}+
\partial_{x_j}\log\theta(p_1+c+x)\frac{\theta_k(p_1)}{\theta_2(p_1)}-
\frac{\theta_{kj}(p_1)}{\theta_2(p_1)}\right.
$$
$$
-\left.
\partial_{x_k}\log\theta(p_2+c+x)\frac{\theta_j(p_2)}{\theta_2(p_2)}-
\partial_{x_j}\log\theta(p_2+c+x)\frac{\theta_k(p_2)}{\theta_2(p_2)}+
\frac{\theta_{kj}(p_2)}{\theta_2(p_2)}\right),
$$

$$
g_{c,c'}^{kj}=\frac{\theta_1(p_1)\theta_1(p_2)}
{\theta_2(p_1)\theta_1(p_2)-\theta_2(p_2)\theta_1(p_1)}
$$
$$
\times\left(
\partial_{x_j}\log\theta(p_1+c+x)\frac{\theta_k(p_1)}{\theta_1(p_1)}+
\partial_{x_k}\log\theta(p_1+c+x)\frac{\theta_j(p_1)}{\theta_1(p_1)}-
\frac{\theta_{kj}(p_1)}{\theta_1(p_1)}
\right.
$$
$$
-\left.
\partial_{x_j}\log\theta(p_2+c+x)\frac{\theta_k(p_2)}{\theta_1(p_2)}-
\partial_{x_k}\log\theta(p_2+c+x)\frac{\theta_j(p_2)}{\theta_1(p_2)}+
\frac{\theta_{kj}(p_2)}{\theta_1(p_2)}\right),
$$

$$
h_{c,c'}^{kj}=\frac{\theta(\Delta+c')}{\theta(\Delta+c'+c+x)}
\bigl(
\partial_{z_k}\partial_{z_j}
-f_{c,c'}^{kj}\partial_{z_1}-g_{c,c'}^{kj}\partial_{z_2}+
2\partial_{z_k}\partial_{z_j}
\log\theta(z)\bigr)
$$
$$
\times
\left(\frac{\theta(z+c+x)}
{\theta(z)}\right)\biggr\vert_{z=\Delta+c'},
$$

$$
[L_{c,c'}(\partial_{z_k}\partial_{z_j}\log\theta(z))]_{12}=
\frac{\theta^2(0)}{\theta(c+c'+x)\theta(c')}
$$
$$
\times\bigl(\partial_{z_k}\partial_{z_j}-
f_{c,c'}^{kj}\partial_{z_1}-g_{c,c'}^{kj}\partial_{z_2}-h_{c,c'}^{kj}
+2\partial_{z_k}\partial_{z_j}
\log\theta(z)\bigr)
\left(\frac{\theta(z+c+x)}{\theta(z)}\right)
\biggr\vert_{z=0},
$$

$$
[L_{c,c'}(\partial_{z_k}\partial_{z_j}
\log\theta(z))]_{21}=
F^{kj}_{c+c',-c'}
\bigl(\alpha_{11}H_{c,c'}^{11}+
\alpha_{12}H_{c,c'}^{12}+
\alpha_{22}H_{c,c'}^{22}+\alpha\bigr),
$$
$$
[L_{c,c'}(\partial_{z_k}\partial_{z_j}
\log\theta(z))]_{22}
=F^{kj}_{c+c',-c'}\bigl(\alpha_{11}F^{11}_{c,c'}+
\alpha_{12}F^{12}_{c,c'}+
\alpha_{22}F^{22}_{c,c'}\bigr)+H^{kj}_{c+c',-c'},
$$
where
$\alpha, \alpha_{kj}\in{\C}$
are such that
$$
\frac{\theta(z-c')\theta(z+c')}{\theta^2(z)}=
\alpha_{11}\partial^2_{z_1}\log\theta(z)+
\alpha_{12}\partial_{z_1}\partial_{z_2}\log\theta(z)+
\alpha_{22}\partial^2_{z_2}\log\theta(z)+\alpha.
$$
}

{\sc Proof of Proposition 2.}
It is easy to verify that the function
$$
 \frac{
\partial_{x_k}\partial_{x_j}\psi
+\partial_{z_k}\partial_{z_j}
\log\theta(z)\psi}
{\exp(-x_1\partial_{z_1}\log\theta(z)
-x_2\partial_{z_2}\log\theta(z))}
$$
has a~ second-order pole on
$\Theta$;
i.e.,
$$
 \partial_{x_k}\partial_{x_j}\psi
+\partial_{z_k}\partial_{z_j}
\log\theta(z)\psi\in M_c(2).
$$
Hence, the operator
$H_{c,c'}^{kj}$
looks like
$$
H_{c,c'}^{kj}=
-\partial_{x_k}\partial_{x_j}+f_{c,c'}^{kj}\partial_{x_1}+
g_{c,c'}^{kj}\partial_{x_2}+
h_{c,c'}^{kj},
$$
and
$F^{kj}_{c,c'}$
is the operator of multiplication by the function
$$
H_{c,c'}^{kj}\psi+F^{kj}_{c,c'}\psi_{c'}=
\partial_{z_k}\partial_{z_j}
\log\theta(z)\psi. \eqno{(7)}
$$
Rewrite ~(7) for
$z\in\Gamma_{c'}$
as follows:
$$
H_{c,c'}^{kj}\psi=
\partial_{z_k}\partial_{z_j}
\log\theta(z)\psi.
$$
By Lemmas~2 and 8, the last equality amounts to
$$
\bigl(\partial_{z_k}\partial_{z_j}
-f_{c,c'}^{kj}\partial_{z_1}-g_{c,c'}^{kj}
\partial_{z_2}-h_{c,c'}^{kj}+
2\partial_{z_k}\partial_{z_j}
\log\theta(z)
\bigr)
\left(\frac{\theta(z+c+x)}{\theta(z)}\right)=0.
\eqno{(8)}
$$
Multiply both sides of~(8) by
$\theta^2(z)$
and put
$z=p_1$ and $z=p_2$
to obtain the following system of linear equations in~
$f_{c,c'}^{kj}$
and
$g_{c,c'}^{kj}$:
$$
f_{c,c'}^{kj}\theta_1(p_1)+
g_{c,c'}^{kj}\theta_2(p_1)=
\frac{
\theta_j(p_1+c+x)\theta_k(p_1)+\theta_k(p_1+c+x)\theta_j(p_1)-
\theta_{kj}(p_1)}
{\theta(p_1+c+x)},
$$
$$
f_{c,c'}^{kj}\theta_1(p_2)+
g_{c,c'}^{kj}\theta_2(p_2)=
\frac{
\theta_j(p_2+c+x)\theta_k(p_2)+\theta_k(p_2+c+x)\theta_j(p_2)-
\theta_{kj}(p_2)}
{\theta(p_2+c+x)}.
$$
Solving it, we find
$f_{c,c'}^{kj}$
and
$g_{c,c'}^{kj}$.
Putting
$z=\Delta+c'$
in~(8)
($\Delta+c'\in\Gamma_{c'}$,
since
$\theta(\Delta)=0$),
we find
$h_{c,c'}^{kj}$.
It follows from Lemmas~2 and 8 that (7) is equivalent to the equality
$$
\bigl(-\partial_{z_k}\partial_{z_j}
+f_{c,c'}^{kj}\partial_{z_1}
+g_{c,c'}^{kj}\partial_{z_2}+h_{c,c'}^{kj}\bigr)
\left(\frac{\theta(z+c+x)}{\theta(z)}\right)
$$
$$
+F_{c,c'}^{kj}
\frac{\theta(z+c+c'+x)\theta(z-c')}{\theta^2(z)}=
2\partial_{z_i}\partial_{z_j}
\log\theta(z)\frac{\theta(z+c+x)}{\theta(z)}.
$$
Putting
$z=0$
in this equality, we obtain
$F_{c,c'}^{kj}$.

We find the remaining entries of the operator
$L_{c,c'}(\partial_{z_k}\partial_{z_j}\log\theta(z))$.
Replace $c'$
with
$-c'$
and
$c$
with
$c+c'$
in ~(7) and multiply both sides by
$\frac{\theta(z-c')}{\theta(z)}$
to obtain
$$
H_{c+c',-c'}^{kj}\psi_{c'}+F^{kj}_{c+c',-c'}
\frac{\theta(z-c')\theta(z+c')}{\theta^2(z)}\psi=
\partial_{z_k}\partial_{z_j}
\log\theta(z)\psi_{c'}.
                            \eqno{(9)}
$$
Let
$$
\frac{\theta(z-c')\theta(z+c')}{\theta^2(z)}=
\alpha_{11}\partial^2_{z_1}\log\theta(z)+
\alpha_{12}\partial_{z_1}\partial_{z_2}\log\theta(z)+
\alpha_{22}\partial^2_{z_2}\log\theta(z)+\alpha,
$$
where
$\alpha, \alpha_{kj}\in{\C}$.
The existence of the constants
$\alpha, \alpha_{kj}$, $k,j=1,2$,
follows from the fact that
the dimension of the space of meromorphic functions with a~pole of order
at most~2 on
$\Theta$
equals~4 (see, for instance, [9]). The functions
$1$
and
$\partial_{z_k}\partial_{z_j}\log\theta(z)$, $k,j=1,2$,
are linearly independent over
${\C}$,
which can be shown, for instance, by means of Nakayashiki's theorem.
The operators
$L_{c,c'}(1)$
and
$L_{c,c'}(\partial_{z_k}\partial_{z_j}\log\theta(z))$, $k,j=1,2$,
are linearly independent, since the leading part of the 11-entries
of these operators, as demonstrated, are equal respectively to
$1$
and
$-\partial_{x_k}\partial_{x_j}$.
Now, (9) amounts to
$$
F^{kj}_{c+c',-c'}(\alpha_{11}H_{c,c'}^{11}+
\alpha_{12}H_{c,c'}^{12}+
\alpha_{22}H_{c,c'}^{22}+\alpha)
\psi
$$
$$
+(F^{kj}_{c+c',-c'}(\alpha_{11}F^{11}_{c,c'}+
\alpha_{12}F^{12}_{c,c'}+
\alpha_{22}F^{22}_{c,c'})+H^{kj}_{c+c',-c'})
\psi_{c'}=
\partial_{z_k}\partial_{z_j}\log\theta(z)\psi_{c'};
$$
whence we obtain
$[L_{c,c'}(\partial_{z_k}\partial_{z_j}\log\theta(z))]_{21}$
and
$[L_{c,c'}(\partial_{z_k}\partial_{z_j}\log\theta(z))]_{22}$,
completing the proof of Proposition~2.

We find the operators
$Z_1$
and
$Z_2$.

{\bf Proposition 3.}
{\sl
The operators
$Z_1$
and
$Z_2$
for the basis
$\psi, \psi_{c'}$
look like
$$
[Z_j]_{11}=-x_1H^{1j}_{c,c'}-x_2H^{j2}_{c,c'}+\partial_{x_j},
\ \ [Z_j]_{12}=-x_1F^{1j}_{c,c'}-x_2F^{j2}_{c,c'},
$$
$$
[Z_j]_{21}=
-x_1[L_{c,c'}(\partial_{z_1}\partial_{z_j}\log\theta(z))]_{21}
$$
$$
-x_2[L_{c,c'}(\partial_{z_j}\partial_{z_2}\log\theta(z))]_{21}+
k_1^j\partial_{x_1}+k_2^j\partial_{x_2}+h^j_{c,c'},
$$
$$
[Z_j]_{22}=
-x_1[L_{c,c'}(\partial_{z_1}\partial_{z_j}\log\theta(z))]_{22}-
x_2[L_{c,c'}(\partial_{z_j}\partial_{z_2}\log\theta(z))]_{22}+
g^j_{c,c'}+2\partial_{x_j},
$$
where
$$
k_1^j=\frac{1}
{\theta_1(p_2)\theta_2(p_1)-\theta_1(p_1)\theta_2(p_2)}
$$
$$
\times
\left(
\frac{\theta_j(p_1-c')\theta(p_1+c+x+c')\theta_2(p_2)}
{\theta(p_1+c+x)}-
\frac{\theta_j(p_2-c')\theta(p_2+c+x+c')\theta_2(p_1)}
{\theta(p_2+c+x)}\right),
$$

$$
k_2^j=\frac{1}
{\theta_1(p_1)\theta_2(p_2)-\theta_1(p_2)\theta_2(p_1)}
$$
$$
\times
\left(
\frac{\theta_j(p_1-c')\theta(p_1+c+x+c')\theta_1(p_2)}
{\theta(p_1+c+x)}-
\frac{\theta_j(p_2-c')\theta(p_2+c+x+c')\theta_1(p_1)}
{\theta(p_2+c+x)}\right),
$$

$$
h^j_{c,c'}(x)=\frac{\theta_j(\Delta)\theta(\Delta+c+x+2c')}
{\theta(\Delta+c')\theta(\Delta+c+x+c')}
$$
$$
-k_1^j\left(
\frac{\theta_1(\Delta+c+x+c')}{\theta(\Delta+c+x+c')}-
\frac{\theta_1(\Delta+c')}{\theta(\Delta+c')}\right)-
k_2^j\left(\frac{\theta_2(\Delta+c'+c+x)}{\theta(\Delta+c'+c+x)}-
\frac{\theta_2(\Delta+c')}{\theta(\Delta+c')}\right),
$$

$$
g_{c,c'}^j(x)=
-\frac{\theta_j(c')}{\theta(c')}+
\frac{\theta_j(c+x+c')}{\theta(c+x+c')}
$$
$$
-\frac{\theta^2(0)}{\theta(c')\theta(c+x+c')}
\bigl(k^j_1\partial_{z_1}+k^j_2\partial_{z_2}+h^j_{c,c'}\bigr)
\left (\frac{\theta(z+c+x)}{\theta(z)}\right)\biggr\vert_{z=0},
$$
$j=1,2.$
}

{\sc Proof of Proposition 3}
We have to find an~operator
$Z_j{\in} \Mat (2,{\cal D})$
such that
$$
[Z_j]_{11}\psi+[Z_j]_{12}\psi_{c'}=
\partial_{z_j}\psi,
\quad
[Z_j]_{21}\psi+[Z_j]_{22}\psi_{c'}=
\partial_{z_j}\psi_{c'},
$$
$$
\partial_{z_j}\psi=\partial_{z_j}
\left(\frac{\theta(z+c+x)}{\theta(z)}\right)
\exp(-x_1
\partial_{z_1}\log\theta(z)
-x_2
\partial_{z_2}\log\theta(z))
$$
$$
-(x_1\partial_{z_1}\partial_{z_j}\log\theta(z)+
x_2\partial_{z_j}\partial_{z_2}\log\theta(z))
\frac{\theta(z+c+x)}{\theta(z)}
\exp(-x_1\partial_{z_1}\log\theta(z)
$$
$$
-x_2\partial_{z_2}\log\theta(z))
=\bigl(-x_1H^{1j}_{c,c'}-x_2H^{j2}_{c,c'}+\partial_{x_j}\bigr)\psi-
(x_1F^{1j}_{c,c'}+x_2F^{j2}_{c,c'})\psi_{c'}.
$$
Therefore,
$$
[Z_j]_{11}=-x_1H^{1j}_{c,c'}-x_2H^{j2}_{c,c'}+\partial_{x_j},
\quad
 [Z_j]_{12}=-x_1F^{1j}_{c,c'}-x_2F^{j2}_{c,c'}.
$$
We find the remaining entries of~
$Z_j$:
$$
\partial_{z_j}\psi_{c'}=
\partial_{z_j}\left(\frac{\theta(z-c')}{\theta(z)}\right)
\frac{\theta(z+c+x+c')}{\theta(z)}
\exp(-x_1
\partial_{z_1}\log\theta(z)
-x_2
\partial_{z_2}\log\theta(z))
$$
$$
+\frac{\theta(z-c')}{\theta(z)}
\partial_{z_j}\left(\frac{\theta(z+c+x+c')}{\theta(z)}\right)
\exp(-x_1
\partial_{z_1}\log\theta(z)-x_2
\partial_{z_2}\log\theta(z))
$$
$$
-\left(x_1\partial_{z_1}\partial_{z_j}\log\theta(z)+
x_2\partial_{z_j}\partial_{z_2}\log\theta(z)\right)
\frac{\theta(z-c')\theta(z+c+x+c')}{\theta^2(z)}
$$
$$
\times\exp(-x_1
\partial_{z_1}\log\theta(z)-x_2
\partial_{z_2}\log\theta(z)).
\eqno{(10)}
$$
Immediate calculations show that
$$
\partial_{z_j}\left(\frac{\theta(z-c')}{\theta(z)}\right)
\frac{\theta(z+c+x+c')}{\theta(z)}
\exp(-x_1
\partial_{z_1}\log\theta(z)-x_2
\partial_{z_2}\log\theta(z))-
\partial_{x_j}\psi_{c'}
$$
$$
=\frac{\theta_j(z-c')\theta(z+c+x+c')-
\theta(z-c')\theta_j(z+c+x+c')}{\theta^2(z)}
$$
$$
\times\exp(-x_1
\partial_{z_1}\log\theta(z)-x_2
\partial_{z_2}\log\theta(z))\in M_c(2);
$$
consequently, there exist functions
$k_1^j, k_2^j, h^j_{c,c'},g^j_{c,c'}\in{\cal O}$
such that
$$
\partial_{z_j}\left(\frac{\theta(z-c')}{\theta(z)}\right)
\frac{\theta(z+c+x+c')}{\theta(z)}
\exp(-x_1
\partial_{z_1}\log\theta(z)-x_2
\partial_{z_2}\log\theta(z))
$$
$$
=\bigl(k_1^j\partial_{x_1}+k_2^j\partial_{x_2}+h^j_{c,c'}\bigr)\psi+
\bigl(g^j_{c,c'}+\partial_{x_j}\bigr)\psi_{c'}.
\eqno{(11)}
$$
To simplify notation, we omit the subscripts
$c$ and $c'$
of
$k_1^j$
and~
$k_2^j$.
We multiply both sides of (11) by
$\theta^2(z)$
and put
$z=p_1$ and $z=p_2$.
We obtain the following system of linear equations in~
$k_1^j$
and
$k_2^j$:
$$
 k_1^j\theta_1(p_1)+k_2^j\theta_2(p_1)=-
\frac{\theta_j(p_1-c')\theta(p_1+c+x+c')}{\theta(p_1+c+x)},
$$
$$
 k_1^j\theta_1(p_2)+k_2^j\theta_2(p_2)=-
\frac{\theta_j(p_2-c')\theta(p_2+c+x+c')}{\theta(p_2+c+x)}.
$$
Solving it, we find
$k_1^j$
and
$k_2^j$.
Putting
$z=\Delta+c'$
in~(11)
($\Delta+c'\in\Gamma_{c'}$),
we find
$h_{c,c'}^j(x)$.
Putting
$z=0$
in~(11), we find
$g_{c,c'}^j(x)$.
From (10) and (11) we obtain
$[Z_j]_{21}$ and $[Z_j]_{22}$.
Proposition~3 is proven.

Let us show that the so-obtained operators generate the ring
$L_{c,c'}(A_\Theta)$.

{\bf Lemma 9.}
{\sl
The operators
$$
L_{\Phi_c}(\partial_{z_k}\partial_{z_j}\log\theta(z)),
\quad
L_{\Phi_c}(\partial_{z_k}\partial_{z_j}\partial_{z_s}\log\theta(z)),
$$
where
$s,j,k=1,2$,
generate the operator ring~
$L_{c,c'}(A_\Theta).$
}

{\sc Proof.}
The dimension of the space of meromorphic functions on~
$X^2$
with pole on
$\Theta$
of multiplicity at most~3 equals~9 (see, for instance, [9]).
Take a~basis for this space which consists of functions of the form
$$
1,\  \frac{f_1}{\theta^3},\dots,\frac{f_8}{\theta^3},
$$
where
$f_i$
is not divisible by
$\theta$,
$i=1,\dots,8$.
As we show below, we can take the basis of the functions
$$
1,\quad  \partial_{z_k}\partial_{z_j}\partial_{z_s}\log\theta(z),
\quad
\partial_{z_1}^3\log\theta(z)+\partial_{z_j}\partial_{z_s}\log\theta(z),
$$
$$
(\partial_{z_1}\partial_{z_2}\log\theta(z))^2-
\partial_{z_1}^2\log\theta(z)\partial_{z_2}^2\log\theta(z),
\eqno{(12)}
$$
where
$k,j,s=1,2$.
By the Lefschetz theorem, the mapping
$(\theta^3:f_1:\dots:f_8)$
defines an~embedding
$F:X^2\rightarrow {\bf CP}^8$
into the projective space. Let
$(y_0:\dots:y_8)$
be homogeneous coordinates in
${\bf CP}^8.$
Then the restrictions of the functions
$\frac{y_1}{y_0},\dots,\frac{y_8}{y_0}$
to
$F(X^2)$
generate the coordinate ring of the affine algebraic variety
$F(X^2\backslash\Theta)$.
Consequently, the functions
$$
\frac{f_1}{\theta^3},\dots,\frac{f_8}{\theta^3}
$$
generate
$A_\Theta$.
We are left with demonstrating that the functions~(12)
are linearly independent over
${\C}$.
This follows from the fact that the operators
$$
L_{c,c'}(1), \
L_{c,c'}(\partial_{z_k}\partial_{z_j}\partial_{z_s}\log\theta(z)),\
L_{c,c'}(\partial_{z_1}^3\log\theta(z)+
\partial_{z_j}\partial_{z_s}\log\theta(z)),
$$
$$
(L_{c,c'}(\partial_{z_1}\partial_{z_2}\log\theta(z))^2-
\partial_{z_1}^2\log\theta(z)\partial_{z_2}^2\log\theta(z)),
$$
where
$k,j,s=1,2$,
are linearly independent, because the leading symbols of the 11-entries
of these operators are equal respectively to
$$
1,\quad  -2\partial_{x_k}\partial_{x_j}\partial_{x_s},
\quad
-2\partial_{x_1}^3-\partial_{x_j}\partial_{x_s},
$$
$$
f_{c,c'}^{22}\partial_{x_1}^3+
\bigl(g_{c,c'}^{22}-2f_{c,c'}^{12}\bigr)\partial_{x_1}^2\partial_{x_2}+
\bigl(f_{c,c'}^{11}-2g_{c,c'}^{12}\bigr)\partial_{x_1}\partial_{x_2}^2+
g_{c,c'}^{11}\partial_{x_2}^3.
$$
Lemma 9 is proven.

{\bf References.}

[1]
 Nakayashiki~A.,
Structure of Baker--Akhiezer modules of
principally polarized Abelian varieties, commuting partial
differential operators and associated integrable systems,
 Duke Math.~J.,
 1991, vol. 62, N. 2, 315--358.

[2]
Dubrovin~B.~A., Krichever~I.~M., and Novikov~S.~P.,
The Schroedinger equation in a~periodic magnetic field and
Riemann surfaces,
Soviet Math. Dokl., 1976, vol. 17, 947--952.

[3]
Veselov~A.~P. and Novikov~S.~P.,
Finite-gap two-dimensional Schr\"odinger operators.
Potential operators,
Soviet Math. Dokl., 1984, vol. 30, 705--708.

[4]
Feldman~J., Knorrer~H., and  Trubowitz~E.,
There is no two-dimensional analogue of Lame's equation,
       Math. Ann., 1992, vol.  294, N. 2, 295--324.

[5]
 Nakayashiki~A.,
      Commuting partial differential operators
and vector  bundles over Abelian varieties,
       Amer.~J. Math., 1994, vol.    116,
       65--100.
[6]
 Mukai~S.,
      Duality between $D(X)$ and $D(\widehat{X})$ with
its application to Picard sheaves, Nagoya Math.~J., 1981,
      vol. 81, 153--175.

[7]
Andreotti~A. and Mayer~A.,
On period relations for Abelian integrals on algebraic curves,
      Ann. Scuola Norm. Sup. Pisa Cl. Sci ~(4), 1967, vol.  21,
       189--238.

[8]
 Krichever~I.~M.,
      The methods of algebraical geometry
      in the theory of nonlinear equations,
       Russian Math. Surveys (1977), vol. 32, no. 6, 185--213.

[9]
 Taimanov~I.~A.,
       Secants of Abelian varieties, theta functions,
      and soliton equations,
      Uspekhi Mat. Nauk, 1997, vol. 52, N. 1, 149--224,

\vskip3mm

{\sl
Sobolev Institute of Mathematics,
Novosibirsk}

mironov@math.nsc.ru

\enddocument